\documentclass{article}
\usepackage{amsfonts,amsmath,amsthm}
\usepackage[latin1]{inputenc}
\newtheorem{proposition}{Proposition}
\newtheorem{prop}{Proposition}
 
\newtheorem{theorem}{Theorem}

\theoremstyle{remark}
\newtheorem{remark}{Remark}
\newtheorem{example}{Example}
\title{Asymptotic results and statistical procedures for time-changed
  L\'evy processes sampled\\ at hitting times}
\author{Mathieu Rosenbaum \quad Peter Tankov \\[0.2cm]  {\normalsize
    Centre de Mathématiques Appliquées,} \\ {\normalsize Ecole Polytechnique, 91128
Palaiseau France} \\ \texttt{\normalsize \{mathieu.rosenbaum,peter.tankov\}@polytechnique.edu}}
\date{}
\begin{document}
\maketitle
\begin{abstract}

\noindent We provide asymptotic results and develop high frequency
statistical procedures for time-changed Lévy processes sampled at
random instants. The sampling times are given by first hitting
times of symmetric barriers whose distance with respect to the
starting point is equal to $\varepsilon$. This setting can be seen
as a first step towards a model for tick-by-tick financial data
allowing for large jumps. For a wide class of Lévy processes, we
introduce a renormalization depending on $\varepsilon$, under
which the Lévy process converges in law to an $\alpha$-stable
process as $\varepsilon$ goes to $0$. The convergence is extended
to moments of hitting times and overshoots. In particular, these
results allow us to construct consistent estimators of the time change
and of the Blumenthal-Getoor index of the underlying Lévy
process. Convergence rates and a central limit theorem are
established under additional assumptions.

\end{abstract}


\section{Introduction}

\noindent In the recent years, a large number of papers has been
devoted to asymptotic results and statistical procedures for
time-changed Lévy processes
\cite{figueroa.lopez.09,figueroa.lopez.10,woerner.07} and more
general semimartingales \cite{aj1,aj2,aj3,aj4,jacod.08}, under
high-frequency discrete sampling.  The classical high frequency
setting consists in observing $n$ values of the process over a
fixed time interval $[0,T]$ at deterministic sampling times
$0=t_0^n<t_1^n<\ldots<t_n^n=T$. Usually, asymptotic results are
given as $n$ goes to infinity and sup$\{t^n_{i+1}-t^n_i\}$ goes to
zero. Motivated by financial applications, many papers focus more
specifically on the asymptotic behavior of volatility estimators.
For example, power variation estimators which are robust to jumps
are studied in \cite{bsw} and \cite{mancini.09}. Since financial data are often seen as
noisy observations of a semimartingale, limit theorems for
volatility estimators under various kinds of perturbations have
also been widely studied, mostly in the case of continuous
semimartingales, see among others
\cite{bhls06,jlmpv07,r07,zma05}.\\

\noindent In this paper we focus on time-changed Lévy models, that
is, we assume that the process of interest $Y$ is given by $Y_t =
X_{S_t}$ where $X$ is a one-dimensional Lévy process and $S$ is a
continuous increasing process (a time change), which plays the
role of the integrated volatility in this setting. Time changed
Lévy models were introduced into financial literature in
\cite{geman} and their estimation from high frequency data with
deterministic sampling was recently addressed in
\cite{figueroa.lopez.09,figueroa.lopez.10}.\\

\noindent In the context of ultra high-frequency financial data,
the assumption of deterministic sampling times is arguably too
restrictive. Several authors have therefore considered volatility
estimation with endogenous sampling times \cite{f2,hjy,li09,rr2}
but so far only in the context of continuous processes.\\

\noindent In this work we assume that the sampling times are given
by first hitting times of symmetric barriers whose distance with
respect to the starting point is equal to $\varepsilon$. More
precisely, the process $Y$ is observed at times
$(T^\varepsilon_i)_{i\geq 0}$ with $T^\varepsilon_0 = 0$ and
$T^\varepsilon_{i+1} = \inf\{t>T^\varepsilon_i : |Y_t -
Y_{T^\varepsilon_i}|\geq \varepsilon\}$ for $i\geq 1$. The
parameter $\varepsilon$ is the parameter driving the asymptotic
and thus we will assume that $\varepsilon$ goes to zero.\\

\noindent This scheme is probably the most simple and common
endogenous sampling scheme. Moreover, in the spirit of \cite{rr2} it can be seen as a first
step towards a model for ultra high frequency financial data
including jump effects.  For example, $Y$ could represent the unobservable
efficient price process and $[-\varepsilon,\varepsilon]$ the bid-ask interval. However, a detailed
financial interpretation of our model is left for further
research. For practical application, the model should in
particular be modified so that the observed values remain on the
tick grid.\\

\noindent Our asymptotic results may more generally open the way
for studying hedging and portfolio strategies with random
endogenous readjustment dates (see e.g. \cite{fukasawa.09,rr3} for
relevant examples in the setting of continuous processes) and for
approximating the solutions of stochastic differential equations
by Euler-type schemes with random discretization dates (see e.g.,
\cite{kohatsu.tankov.09,rubenthaler}).\\

\noindent We focus on the class of Lévy processes such that for a
suitable $\alpha$, the rescaled process $(X^\varepsilon_t)_{t\geq
0} := (\varepsilon^{-1} X_{\varepsilon^\alpha t})_{t\geq 0}$
converges in law to an $\alpha$-stable Lévy process $X^*$ as
$\varepsilon$ goes to zero. This class turns out to be rather
large, and contains in particular all Lévy processes with
non-zero diffusion component, all finite variation Lévy processes
with non-zero drift and also most parametric Lévy models found in
the literature. We show that for such Lévy processes the
moments of first exit times from intervals, and certain
functionals of the overshoot converge to the corresponding
functionals of the limiting stable process, which are often known
explicitly.\\

\noindent These findings, which are of interest in their own
right, allow us to prove limit theorems for quantities of the form
$$V^\varepsilon(f)_t=\sum_{T^\varepsilon_i \leq t}
f\big(\varepsilon^{-1}(Y_{T^\varepsilon_i}-Y_{T^\varepsilon_{i-1}})\big),$$
leading to consistent estimators of the time change and of the
characteristics of $X$ that are preserved by the limiting
procedure, such as, for example, the Blumenthal-Getoor index of
jump activity.  In some cases, we are able to quantify the rate of
convergence of the functionals of the rescaled process
$X^\varepsilon$ to the corresponding functionals of the limiting
stable process $X^*$.  From this, convergence rates and central
limit theorems for our estimators can be deduced.\\

\noindent The paper is organized as follows. In Section
\ref{conv}, we study the convergence in law of the properly
rescaled underlying L\'evy process $X$ as $\varepsilon$ goes to
zero. Asymptotic results for the first exit time and the overshoot
(more precisely we study the value of the process at the first
exit time which is directly related to the overshoot) are given in
Section \ref{rate}. The law of large numbers for
$V^\varepsilon(f)$ is stated  in Section \ref{lln}, where we also
discuss statistical applications.  Finally, a multidimensional central
limit theorem is given in Section \ref{secclt}. The proofs are relegated to
Section \ref{proofs}.

\section{Convergence of the rescaled process}\label{conv}

\noindent In this section, we give results on the convergence in
law of the properly rescaled process $X$ as $\varepsilon$ goes to
zero. The convergences in law are given in the Skorohod space, for
the usual Skorohod topology. These results will be essential for
proving the law of large numbers and the central limit theorem.
Let us first recall the definition of a strictly stable process and introduce other useful notation.

\paragraph{Preliminaries and notation}
We denote by $(A,\nu,\gamma)$ the characteristic triplet
of the one-dimensional Lévy process $X$, with respect to a
truncation function $h$. This means that via the Lévy-Khintchine formula, the characteristic function of $X_t$ is 
$$
E[e^{iuX_t}] = e^{t\psi(u)},\quad \psi(u) = -\frac{Au^2}{2} + i\gamma u + \int_{\mathbb R}(e^{iux} - 1 - iu h(x))\nu(dx).
$$
Unless otherwise specified, we assume
$h(x) = -1 \vee x \wedge 1$. \\

\noindent A Lévy process $X$ is called strictly $\alpha$-stable for $\alpha\in (0,2]$ if $X_t$ has a strictly $\alpha$-stable distribution for all $t$. This happens if and only if $X$ is selfsimilar, that is, 
$$
\forall a>0,\quad \left(\frac{X_{at}}{a^{1/\alpha}}\right)_{t\geq 0}
{=}\, (X_t)_{t\geq 0},\ \text{in law.}
$$
As recalled in the following proposition, strictly stable Lévy processes can be described in terms of their characteristic triplet. 
\begin{prop}[Theorems 14.3, 14.7 in \cite{sato}] Let $X$ be a Lévy process with characteristic triplet $(A,\nu,\gamma)$. 
\begin{enumerate}
\item $X$ is strictly $2$-stable if and only if $\nu = 0$ and $\gamma=0$. 
\item $X$ is strictly $\alpha$-stable with $1<\alpha<2$ if and only if $A=0$, $\nu$ has a density of the form
\begin{align}
\nu(x) = \frac{c_+}{|x|^{1+\alpha}}1_{x>0} + \frac{c_-}{|x|^{1+\alpha}}1_{x<0},\label{stabledens.eq}
\end{align}
and $\gamma_c=0$ where $\gamma_c:=\gamma - \int_{\mathbb R} (h(x)-x) \nu(dx)$ is the third component of the characteristic triplet of $X$ with respect to the truncation function $h(x)=x$.
\item $X$ is strictly $1$-stable if and only $A=0$ and $\nu$ has a density of the form
$$\nu(x) = \frac{c}{|x|^{2}}.$$
\item $X$ is strictly $\alpha$-stable with $0<\alpha<1$ if and only if
  $A=0$, $\nu$ has a density of the form \eqref{stabledens.eq} and
  $\gamma_0=0$, where $\gamma_0:=\gamma - \int_{\mathbb R} h(x)
  \nu(dx)$ is the third component of the characteristic triplet of $X$
  with respect to the truncation function $h = 0$. 
\end{enumerate}
\end{prop}

\noindent  For $\alpha \in (0,2]$ and
$\varepsilon>0$, we define the rescaled Lévy process
$X^\varepsilon$ via $X^\varepsilon_t :=
\varepsilon^{-1}X_{\varepsilon^\alpha t}$, $t\geq 0$. The first
exit time by the rescaled process from the interval $(-1,1)$ will
be denoted by $\tau^\varepsilon_1:=\inf\{t\geq 0:
|X^\varepsilon_t|\geq 1\}$. This time is directly related to the
first exit time by the original process from the interval
$(-\varepsilon,\varepsilon)$:
$$
\inf\{t\geq 0: |X_t|\geq \varepsilon\} = \varepsilon^{\alpha}
\tau^\varepsilon_1.
$$
Similarly, $X^\varepsilon_{\tau^\varepsilon_1}$ is equal to
$\varepsilon^{-1}$ times the value of $X$ at first exit from
$(-\varepsilon,\varepsilon)$. From the Lévy-Khintchine formula it
is easy to see that the characteristic triplet
$(A^\varepsilon,\nu^\varepsilon,\gamma^\varepsilon)$ of
$X^\varepsilon$ is given by
\begin{align}
A^\varepsilon &= A \varepsilon^{\alpha - 2};\label{aeps}\\
\nu^\varepsilon(B) &= \varepsilon^\alpha \nu(\{x:x/\varepsilon \in B\}),\quad B\in \mathcal B(\mathbb R);\label{nueps}\\
\gamma^\varepsilon &= \varepsilon^{\alpha-1} \big\{\gamma +
\int_{\mathbb R} \nu(dx) (\varepsilon
h(x/\varepsilon)-h(x))\big\}.\label{gammaeps}
\end{align}


\paragraph{Assumptions} To be able to prove the convergence of the properly rescaled
process, we introduce two assumptions on the Lévy measure which
will sometimes be imposed in the sequel:
\begin{itemize}
\item[($\mathbf{H}$-$\alpha$)] The Lévy measure $\nu$ has a density $\nu(x) = \frac{g(x)}{|x|^{1+\alpha}}$, where $g$ is a nonnegative measurable function admitting left and right limits at zero:
$$
c^+:=\lim_{x\downarrow 0} g(x),\quad c^-:=\lim_{x\uparrow 0} g(x),
$$
with $c_+ + c_->0$.
\item[($\mathbf{H'}$-$\alpha$)] The Lévy measure $\nu$ satisfies ($\mathbf{H}$-$\alpha$) and additionnally $c_+c_->0$ and the function $g$ is left- and right-Hölder continuous at zero with exponent $\theta>\alpha/2$:
$$
\limsup_{x\downarrow 0}\frac{|g(x)-c_+|}{|x|^\theta} < \infty\quad
\text{and}\quad \limsup_{x\uparrow 0}\frac{|g(x)-c_-|}{|x|^\theta}
< \infty.
$$
\end{itemize}

\paragraph{Convergence in law of the rescaled process} We now establish a set of alternative conditions under which the
rescaled process $X^\varepsilon$ converges in law to a strictly
stable process as $\varepsilon \to 0$. In the sequel, we will
always work under one of these alternative assumptions. The
following proposition, therefore, also serves as the definition of
the limiting process $X^*$ and of the scaling parameter $\alpha$
depending on the characteristics of $X$.
\begin{proposition}${}$\label{cvginlaw.prop}
\begin{enumerate}
\item \label{brownian.cond} Let $\alpha=2$ and $A>0$. Then the process $X^\varepsilon$ converges in law to a Lévy process $X^*$ with characteristic triplet $(A,0,0)$, that is, to a Brownian motion with variance $A$ at time $t=1$.
\item \label{drift.cond} Let $\alpha=1$ and assume that $X$ has finite variation (that is, $A=0$ and $\int_{|x|\leq 1} |x|\nu(dx)<\infty$) and nonzero drift: $\gamma_0:=\gamma - \int_{\mathbb R} h(x) \nu(dx) \neq 0$. Then the process $X^\varepsilon$ converges in law to the (deterministic) Lévy process $X^*$ with characteristic triplet $(0,0,\gamma_0)$.
\item \label{infvar.cond} Let $1<\alpha<2$ and assume that $A=0$ and that the Lévy measure $\nu$ satisfies the condition ($\mathbf{H}$-$\alpha$). Then the process $X^\varepsilon$ converges in law to a strictly $\alpha$-stable Lévy process $X^*$ with Lévy density
\begin{align}
\nu^*(x) = \frac{c_+ 1_{x>0} + c_-1_{x<0}}{|x|^{1+\alpha}}.\label{stabledens}
\end{align}
\item \label{cauchy.cond} Let $\alpha=1$ and assume that $A=0$ and that the Lévy measure $\nu$ satisfies the condition ($\mathbf{H}$-$\alpha$) with $c^+=c^-:=c$ and with the function $g$ satisfying
$$
\int_0^1 \frac{|g(x)-g(-x)|dx}{x}<\infty.
$$
Then the process $X^\varepsilon$ converges in law to a Lévy process $X^*$ with characteristic triplet $(0,\nu^*,\gamma^*)$, where $\gamma^* = \gamma - \int_0^\infty \frac{g(x)-g(-x)}{x^2}h(x)dx$ and $\nu^*$ has Lévy density
$$
\nu^*(x) = \frac{c}{|x|^2},
$$
that is, to a strictly $1$-stable Lévy process.
\item \label{finvar.cond} Let $0<\alpha<1$ and assume that $A=0$, the process has zero drift: $\gamma - \int_{\mathbb R} h(x) \nu(dx) = 0$  and that the Lévy measure $\nu$ satisfies the condition ($\mathbf{H}$-$\alpha$). Then the process $X^\varepsilon$ converges in law to a strictly $\alpha$-stable Lévy process $X^*$ with Lévy density \eqref{stabledens}.
\end{enumerate}
\end{proposition}

\begin{remark}
This result is closely related to the convergence of tempered stable
processes to stable processes studied in \cite{rosinski.04}. More
precisely, in Theorem 3.1 of \cite{rosinski.04}, Rosi\'nski proves the
results of parts \ref{infvar.cond}, \ref{cauchy.cond} and
\ref{finvar.cond} under the additional assumption that the function
$g$ is completely monotone (but in the multidimensional setting). 
\end{remark}

\begin{remark}
The different alternative cases contain the main parametric models found in finance literature. We list several examples below.
\begin{itemize}
\item All models with a nonzero diffusion component (e.g., the models of Merton \cite{merton} and Kou \cite{kou}) satisfy Condition \ref{brownian.cond}.
\item The variance gamma model \cite{madan98} with nonzero drift satisfies Condition \ref{drift.cond}.
\item The normal inverse gaussian process (NIG), see \cite{bns_nig}, satisfies Condition \ref{cauchy.cond}. This can be seen directly from the form of the Lévy density
$$
\nu(x) = \frac{C}{|x|}e^{Ax} K_1(B|x|),
$$
where $A$, $B$ and $C$ are constants and $K_1$ is the modified Bessel function of the second kind, which satisfies $K_1(x) \sim \frac{1}{x}$ for $x\downarrow 0$.
\item The CGMY process, see \cite{finestructure}, that is, a Lévy process with no diffusion component and a Lévy density of the form
\begin{align}
\nu(x) = \frac{C e^{-\lambda_-|x|}}{|x|^{1+\alpha}}1_{x<0} + \frac{C e^{-\lambda_+|x|}}{|x|^{1+\alpha}}1_{x>0}\label{cgmy.eq},
\end{align}
satisfies Condition \ref{infvar.cond} if $1<\alpha<2$, Condition \ref{cauchy.cond} if $\alpha=1$, Condition \ref{drift.cond} if $\alpha<1$ and the process has nonzero drift, and Condition \ref{finvar.cond} if $\alpha<1$ and the drift is zero.
\end{itemize}
\end{remark}

\section{Asymptotic results for the first exit time and the overshoot of Lévy processes out of small
intervals}\label{rate}

In this section, our aim is to study the
first exit time and the overshoot corresponding to the exit of $X$
from the interval $(-\varepsilon,\varepsilon)$. In order to work
with quantities of order 1, we formulate our results in terms of
$\tau^\varepsilon_1$ and $X^\varepsilon_{\tau^\varepsilon_1}$.

\paragraph{Convergence for the first exit time and overshoot}
We define $\tau^*$ as the first exit time by the limiting
process $X^*$ from the interval $(-1,1)$. Observe that $\tau^*$
admits moments of any order. When $X^*$ is a nontrivial
$\alpha$-stable process with $0<\alpha<2$, $\tau^*$ is dominated
by the time of the first jump of $X^*$ greater than $2$ in
absolute value, which has exponential distribution. In the case
$\alpha=2$ (Brownian motion) this is a classical result, see for
example \cite{ds,d}.

\begin{proposition}${}$\label{cvgtimes.prop}
Let $X$ be a Lévy process satisfying one of the conditions 1--5 of
Proposition \ref{cvginlaw.prop} and let $f$ be a bounded continuous
function on $\mathbb{R}$. Then
\begin{enumerate}
\item  $(\tau^\varepsilon_1,X^\varepsilon_{\tau^\varepsilon_1})$ converges in law to $(\tau^*_1,X^*_{\tau^*_1})$ as $\varepsilon \downarrow 0$.
\item  $\lim_{\varepsilon \downarrow 0}E[(\tau^\varepsilon_1)^kf(X^\varepsilon_{\tau^\varepsilon_1})] = E[(\tau^*_1)^kf(X^*_{\tau^*_1})]$ for all $k\geq 1$.
\end{enumerate}
\end{proposition}

\begin{remark}
The weak convergence of the $X^\varepsilon_{\tau^\varepsilon_1}$ under
Conditions \ref{brownian.cond} or \ref{drift.cond} of Proposition
\ref{cvginlaw.prop} (actually in these two cases
$|X^\varepsilon_{\tau^\varepsilon_1}|\to 1$) is a known result
\cite{doney.maller.02}. See also \cite[Theorem 5.16]{kyprianou} for a related result
in the context of subordinators.
\end{remark}

\begin{remark}
The moments of the exit time and the law of the overshoot for the
limiting strictly stable process are often known explicitly. 
\begin{itemize}
\item Under Condition \ref{brownian.cond} of Proposition
\ref{cvginlaw.prop}, the limiting process is a Brownian motion, so
$X^*_{\tau^*_1}$ equals $1$ or $-1$ with probability $\frac{1}{2}$ and
the law of $\tau^*_1$ is well known (see e.g., exercise II.3.10 in \cite{revuz.yor.99}). 
\item Under
Condition \ref{drift.cond} of Proposition
\ref{cvginlaw.prop}, the limiting process is deterministic, so
$\tau^*_1 = \frac{1}{|\gamma_0|}$ and $X^*_{\tau^*_1} =
\text{sgn}\,\gamma_0$. 
\item Under Conditions \ref{infvar.cond}--\ref{finvar.cond}, the first and
second moments of the hitting time $\tau^*_1$ are given in
\cite{getoor.61} for the symmetric case, and the law of the overshoot
is computed in \cite{blumenthal.al.61} for the symmetric case and in
\cite{rogozin.72} for the general case. 
\end{itemize}
\end{remark}

\paragraph{Rates of convergence for the first exit times and overshoots}

We now compute the rates of convergence of
$E[\tau^\varepsilon_1]$ to $E[\tau^*_1]$ 	and of $E[f(X^\varepsilon_{\tau^\varepsilon_1})]$ to $E[f(X^*_{\tau^*_1})]$. These results either guarantee the asymptotic
normality of the estimators provided in Section \ref{lln} or allow to establish a convergence rate or an error bound for these estimators
in the cases when the bias asymptotically dominates the variance.

\begin{proposition}\label{rate_exit_times} ${}$
\begin{enumerate}
\item\label{clt1} Let $X$ be a Lévy process satisfying Condition
  \ref{brownian.cond} of Proposition \ref{cvginlaw.prop} such that its
  Lévy measure $\nu$ satisfies $\int_{|x|\leq 1} |x|\nu(dx)<\infty$ and let $f$ be a bounded Lipschitz function on $\mathbb R$ with $f(-1)=f(1)$.  Then
$$
\lim_{\varepsilon \downarrow 0}\varepsilon^{-1}(E[\tau^\varepsilon_1] - E[\tau^*_1])=0 \quad \text{and}\quad \lim_{\varepsilon \downarrow 0}\varepsilon^{-1}(E[f(X^\varepsilon_{\tau^\varepsilon_1})] - E[f(X^*_{\tau^*_1})])=0.
$$
\item \label{ratecond}
Let $X$ be a Lévy process satisfying Condition \ref{drift.cond} of
Proposition \ref{cvginlaw.prop} such that its Lévy measure $\nu$
satisfies $\int_{|x|\leq 1}|x|^\beta \nu(dx)$ for some $\beta \in(0,1)$, and let $f$ be a bounded Lipschitz function on $\mathbb R$. Then
\begin{align}
\lim_{\varepsilon \downarrow 0}\varepsilon^{-(1-\beta-\delta)}(E[\tau^\varepsilon_1] - E[\tau^*_1])=0\label{beta1}
\end{align}
and 
\begin{align}
\lim_{\varepsilon \downarrow 0}\varepsilon^{-(1-\beta-\delta)}(E[f(X^\varepsilon_{\tau^\varepsilon_1})] - E[f(X^*_{\tau^*_1})])=0\label{beta2}
\end{align}
for all $\delta>0$.
\item\label{clt2} Let $X$ be a Lévy process satisfying Condition
  \ref{infvar.cond} of Proposition \ref{cvginlaw.prop}, Assumption
  $(\mathbf H'$-$\alpha)$ and the condition
\begin{align}
\gamma= \int_{\mathbb R}\left(h(x)\frac{d\nu}{d\nu^*}(x)-x\right)\nu^*(dx).\label{driftconst}
\end{align}

\centerline{or}

let $X$ be a Lévy process satisfying Condition \ref{cauchy.cond} of Proposition \ref{cvginlaw.prop} and Assumption $(\mathbf H'$-$\alpha)$

\centerline{or}

let $X$ be a Lévy process satisfying Condition \ref{finvar.cond} of Proposition \ref{cvginlaw.prop} and Assumption $(\mathbf H'$-$\alpha)$. Let $f$ be a bounded continuous function on $\mathbb R$. 

Then
\begin{align}
\lim_{\varepsilon \downarrow 0}\varepsilon^{-\alpha/2}(E[\tau^\varepsilon_1] - E[\tau^*_1])=0
 \quad \text{and}\quad \lim_{\varepsilon \downarrow 0}\varepsilon^{-\alpha/2}(E[f(X^\varepsilon_{\tau^\varepsilon_1})] - E[f(X^*_{\tau^*_1})])=0.
\label{ratetau.eq}
\end{align}
\end{enumerate}
\end{proposition}

\begin{remark}
As we shall see below, Conditions \ref{clt1} and \ref{clt2} lead to a
central limit theorem for the estimators constructed in the following
sections, while Condition \ref{ratecond} provides a convergence rate
without ensuring asymptotic normality.  A natural question is what
happens in the case where the Lévy process satisfies Condition
\ref{infvar.cond} of Proposition \ref{cvginlaw.prop} but the drift
constraint \eqref{driftconst} is not satisfied. In this case, we have been unable to obtain a convergence rate, due to unsufficient regularity of the functions of type $E^x[\tau^*_1]$ and $E^x[f(X^*_{\tau^*_1})]$. However the following example shows that the estimate \eqref{ratetau.eq} may not hold in this case, and therefore one cannot hope to obtain a limit theorem without bias. \\

\noindent Let $X$ be a Lévy process with characteristic triplet
$(0,\nu,\gamma_c)$ with respect to the truncation function $h(x)=x$
and $\nu$ given by \eqref{stabledens.eq} with $c_+=c_-$ and
$1<\alpha<2$. Assume $\gamma_c>0$ (hence the drift constraint is not
satisfied) and let $f(x) = 1_{(1,\infty)}(x) + x 1_{(0,1]}(x)$. The process $X^*$ then has the characteristic triplet $(0,\nu,0)$ (with respect to the same truncation function), and the function $u(x):=E^x[f(X^*_{\tau^*_1})]$ is given by (see \cite{blumenthal.al.61}),
$$
u(x) = 2^{1-\alpha}\Gamma(\alpha)\left[\Gamma\left(\frac{\alpha}{2}\right)\right]^{-2} \int_{-1}^x (1-u^2)^{\alpha/2-1}du
$$
for $|x|<1$ and $u(x) = f(x)$ for $|x|\geq 1$. Observe that for $|x|<1$, 
$$
u'(x)\geq 2^{1-\alpha}\Gamma(\alpha)\left[\Gamma\left(\frac{\alpha}{2}\right)\right]^{-2}:= C
$$
and (this is shown in \cite{blumenthal.al.61}) 
$$
\int_{\mathbb R}\left\{u(x+z)-u(x)-zu'(x)\right\}\nu(dz) = 0.
$$
Using this identity in the Itô formula applied to $u(X^\varepsilon_t)$
between $t=0$ and $t=\tau^\varepsilon_{\delta}$ for $\delta\in(0,1)$ (to avoid regularity issues), and taking the expectation, we get 
$$
E[u(X^\varepsilon_{\tau^\varepsilon_{\delta}}) - u(0)] = \varepsilon^{\alpha-1}\gamma_c E\left[\int_0^{\tau^\varepsilon_{\delta}} u'(X^\varepsilon_s)\right] \geq C \varepsilon^{\alpha-1}\gamma_c E[\tau^\varepsilon_{\delta}],
$$
which is equivalent to
$$
E[u(\delta X^{\varepsilon\delta}_{\tau^{\varepsilon\delta}_1}) - u(0)] \geq C
\varepsilon^{\alpha-1} \delta^\alpha \gamma_c E[\tau^{\varepsilon\delta}_1].
$$
With the notation $\rho = \varepsilon\delta$, this gives
$$
E[u(\delta X^{\rho}_{\tau^{\rho}_1}) - u(0)] \geq C
\rho^{\alpha-1} \delta \gamma_c E[\tau^{\rho}_1].
$$
Taking the limit $\delta\to 1$ then yields
$$
E[f(X^\rho_{\tau^\rho_1})-f(X^*_{\tau^*_1})] = E[u(X^\rho_{\tau^\rho_{1}}) - u(0)] \geq C \rho^{\alpha-1}\gamma_c E[\tau^\rho_1],
$$
which is bounded from below by $\rho^{\alpha-1}$ times a positive constant since $E[\tau^\rho_1]$ converges to $E[\tau^*_1]$.
\end{remark}


\section{Law of large numbers and statistical applications}\label{lln}
In this section we give the law of large numbers for the the processes
of the form
$$V^{\varepsilon}(f)_t=\sum_{T_i^{\varepsilon}\leq t}f\big(\varepsilon^{-1}(Y_{T^{\varepsilon}_{i}}-Y_{T^{\varepsilon}_{i-1}})\big),$$
where $f$ is a bounded continuous function on $\mathbb{R}$.
Let $$m(f)=\frac{E[f(X^*_{\tau^*_1})]}{E[\tau^*_1]}.$$ 
\begin{theorem}\label{overshootucp.prop}
Let $X$ be a Lévy process with characteristic triplet
$(A,\nu,\gamma)$, satisfying one of the conditions 1--5 of
Proposition \ref{cvginlaw.prop}. Let $f$ be a bounded continuous function on $\mathbb{R}$. Then
\begin{align}
\lim_{\varepsilon\downarrow 0}\varepsilon^\alpha
V^{\varepsilon}(f)_t = m(f)S_t\label{overshootucp.eq}
\end{align}
in probability, uniformly on compact sets in $t$ (ucp).
\end{theorem}

\noindent As shown in the following examples, this result can be
in particular used to build estimators of relevant quantities such
as the time change or the Blumenthal-Getoor index.

\begin{example}[Estimation of the time change]
Assume that the parameters of the underlying Lévy process are known. In our model, the time change can be recovered
simply from the times $(T^\varepsilon_i)$ as $\varepsilon\to 0$, by taking $f=1$, which gives,
\begin{align}\label{poissonucp.eq}
S_t = \lim_{\varepsilon\downarrow 0}\varepsilon^\alpha
V^\varepsilon(1)_t E[\tau^*_1].
\end{align}
\end{example}

\begin{example}[Estimation of the Blumenthal-Getoor index for the time-changed CGMY process]
Let $X$ be the CGMY process \eqref{cgmy.eq} with $1<\alpha<2$.
Including the constant $C$ into the time change, we can assume
$C=1$ with no loss of generality. In this case, the limiting
process $X^*$ is a symmetric $\alpha$-stable process and has Lévy
density $\nu^*(x) = \frac{1}{|x|^{1+\alpha}}$. Our method allows
therefore to estimate the Blumenthal-Getoor index $\alpha$ of the
process $X$. The coefficients $\lambda_+$ and $\lambda_-$ cannot
be identified from the trajectory of the process over a finite
time interval, even in the case of continuous observation.\\

\noindent The law of the symmetric stable process at the first
exit time from an interval is well known in the literature
\cite{blumenthal.al.61,getoor.61}: $X^*_{\tau^*_1}$ has density
$$
\mu(y) =
\frac{1}{\pi}\sin\big(\frac{\pi\alpha}{2}\big)|y|^{-1}(y^2-1)^{-\frac{\alpha}{2}},\quad
|y|\geq 1.
$$
and
\begin{align}
E[\tau^*_1] =
\frac{\sqrt{\pi}}{2^{\alpha}\Gamma\big(1+\frac{\alpha}{2}\big)}.\label{etau}
\end{align}
With $f(x) = \frac{1}{|x|^\beta}\wedge 1$, $\beta\geq 0$ we easily get
$$
E[f(X^*_{\tau^*_1})] = \int_{|y|\geq 1} \frac{\mu(y)}{|y|^\beta}dy
= \frac{\Gamma\big(\frac{\alpha}{2} +
\frac{\beta}{2}\big)}{\Gamma\big(\frac{\alpha}{2}
\big)\Gamma\big(1 + \frac{\beta}{2}\big)},
$$
where $\Gamma$ is the gamma function,
and in particular for $\beta=2$, $E[(X^*_{\tau^*_1})^{-2}] =
\frac{\alpha}{2}$. Combining \eqref{overshootucp.eq} and
\eqref{poissonucp.eq}, we then obtain a consistent estimator of
$\alpha$:
$$
\alpha = 2 \lim_{\varepsilon \downarrow 0}
\frac{V^\varepsilon(f)_t}{V^\varepsilon(1)_t},\quad f(x) =
\frac{1}{x^2}\wedge 1.
$$

\end{example}

\section{Central limit theorem and convergence rates for estimators}\label{secclt}

We now turn to the central limit theorem. The following result
establishes the rate of convergence and asymptotic normality of
the renormalized error in \eqref{overshootucp.eq}.
\begin{theorem}\label{clt}
Assume that the
time-change $S$ defining $Y$ is independent of the underlying
L\'evy process $X$.

\medskip

\noindent Let $X$ be a Lévy process satisfying Condition \ref{clt1}  of
Proposition \ref{rate_exit_times} and  let $d\in\mathbb{N}^{*}$ and $f_1,\ldots,f_d$
be bounded Lipschitz functions on $\mathbb{R}$ satisfying $f_i(1)=f_i(-1)$ for $i=1,\dots, d$

\centerline{or}

\noindent let $X$ be a Lévy process satisfying Condition \ref{clt2}  of
Proposition \ref{rate_exit_times} and  let $d\in\mathbb{N}^{*}$ and $f_1,\ldots,f_d$
be bounded continuous functions on $\mathbb{R}$.

\medskip

\noindent Define
$R^\varepsilon_t =
(R^\varepsilon_{t,1},\ldots,R^\varepsilon_{t,d})$ with
$$
R^\varepsilon_{t,j}= \varepsilon^{-\alpha/2}(\varepsilon^\alpha
V^{\varepsilon}(f_j)_t-m(f_j)S_t).$$ Then, as $\varepsilon$ goes
to zero, $R^\varepsilon$ converges in law to $B\circ S$, for the
usual Skorohod topology, with $B$ a continuous centered
$\mathbb{R}^d-$valued Gaussian process with independent
increments, independent of $S$, such that
$E[B_{t,j}B_{t,k}]=(t/(E[\tau^*_1])C_{j,k}$ with
$$C_{j,k}=\emph{Cov}[f_j(X^{*}_{\tau_1^{*}})-m(f_j)\tau^*_1,f_k(X^{*}_{\tau_1^{*}})-m(f_k)\tau^*_1].$$
\end{theorem}

\medskip

\noindent Under Condition \ref{ratecond} of Proposition \ref{rate_exit_times}, $\tau^*_1$ et $X^*_{\tau^*_1}$ are deterministic, and therefore a central limit theorem cannot be established. In this case, we can only provide an upper bound on the error of the estimators. 
\begin{proposition}\label{noclt}
Let $X$ be a Lévy process satisfying Condition \ref{ratecond} of Proposition \ref{rate_exit_times}, and let $f$ be a real bounded Lipschitz function  on $\mathbb R$. Then, for every $\delta>0$,
$$
\varepsilon^{-(1-\beta-\delta)\vee -\frac{1}{2}} \{\varepsilon V^\varepsilon(f)_t - m(f)S_t\} \to 0
$$
as $\varepsilon\to 0$, in probability uniformly in $t$ on compacts. 
\end{proposition}
\section{Proofs}\label{proofs}

We give in this section the proofs of the preceding results.

\subsection{Proof of Proposition \ref{cvginlaw.prop}}

Let $(A^*,\nu^*,\gamma^*)$ denote the characteristic triplet of
the limiting process. By corollary VII.3.6 in
\cite{jacodshiryaev}, in order to prove the convergence in law, we
need to check that
\begin{align}
\gamma^\varepsilon &\to \gamma^*;\label{gammacond}\\
A^\varepsilon + \int_{\mathbb R} h^2(x)\nu^{\varepsilon}(dx) &\to A^* + \int_{\mathbb R} h^2(x)\nu^{*}(dx);\label{acond}\\
\text{and}\quad \int_{\mathbb R} f(x)\nu^\varepsilon(dx) &\to
\int_{\mathbb R} f(x)\nu^*(dx)\label{nucond}
\end{align}
for every continuous bounded function $f$ which is zero in a
neighborhood of zero.

\paragraph{Part 1} We first check \eqref{gammacond}. Using the
explicit form of the truncation function, we get, for $\varepsilon<1$,
$$
|\gamma^\varepsilon| \leq \varepsilon|\gamma| +
\varepsilon\int_{|x|>1}\nu(dx) + \varepsilon
\int_{\varepsilon<x\leq 1} (x-\varepsilon)\nu(dx)+ \varepsilon
\int_{-1\leq x<-\varepsilon} (\varepsilon-x)\nu(dx).
$$
The convergence of the first two terms to zero is evident; for the
third term it is the consequence of the dominated convergence
theorem, because the integrand $\varepsilon
(x-\varepsilon)1_{\varepsilon<x\leq 1}$ converges to zero and is
bounded from above by $x^2 1_{0<x \leq 1}$, and the fourth term is
treated similarly to the third one. Therefore, $\gamma^\varepsilon
\to 0 = \gamma^*$.\\

\noindent To prove \eqref{acond}, we observe that $A^\varepsilon
\to A$ and moreover
$$
\int_{\mathbb R} h^2(x)\nu^\varepsilon(dx) = \int_{|x|\leq
\varepsilon} x^2 \nu(dx) + \varepsilon^2 \int_{|x|>1}\nu(dx) +
\varepsilon^2 \int_{\varepsilon<|x|\leq 1} \nu(dx).
$$
For the first two terms the convergence to zero is evident, and
for the last one we can once again apply the dominated convergence
theorem using the fact that $\varepsilon^2 1_{\varepsilon<|x|\leq
1} \leq x^2 1_{0<|x|\leq 1}$.\\

\noindent For the condition \eqref{nucond}, assume $f(x)=0$ for
$|x|\leq \delta$. Then we can again decompose
$$
\int_{\mathbb R}f(x)\nu^\varepsilon(dx) = \varepsilon^2 \int_{
\delta \varepsilon<|x|\leq 1}f(x/\varepsilon)\nu(dx) +
\varepsilon^2 \int_{|x|>1}f(x/\varepsilon)\nu(dx),
$$
and apply the dominated convergence theorem to the first term, to
show that the limit is zero.

\paragraph{Part 2} The proof of this part is a minor modification of part 1, so we omit it to save space.

\paragraph{Conditions \eqref{acond} and \eqref{nucond} in parts 3, 4 and 5} To prove \eqref{acond}, we fix $\eta>0$ such that $g(x)$ is bounded on $[-\eta,\eta]$. Then
\begin{align*}
\lim_{\varepsilon\downarrow 0}\int_{\mathbb R} h^2 (x)
\nu^\varepsilon(dx) &= \lim_{\varepsilon\downarrow 0}
\varepsilon^\alpha \int_{|x|\leq \eta}
\frac{h^2(x/\varepsilon)g(x)dx}{|x|^{1+\alpha}}\\ &=
\lim_{\varepsilon\downarrow 0} \int_{|x|\leq \eta/\varepsilon}
\frac{h^2(x)g(\varepsilon x)dx}{|x|^{1+\alpha}} = \int_{\mathbb R}
h^2(x)\nu^*(dx),
\end{align*}
where in the last equality we use the dominated convergence theorem. The
condition \eqref{nucond} is shown in a similar manner.

\paragraph{Condition \eqref{gammacond} in part 3} Since $\alpha >1$ and $h$ is bounded, for every $\eta>0$,
$$
\lim_{\varepsilon\downarrow 0} \gamma^\varepsilon =
\lim_{\varepsilon\downarrow 0} \varepsilon^{\alpha-1}
\int_{|x|\leq \eta} \nu(dx) (\varepsilon h(x/\varepsilon)-h(x))
$$
Since $g$ has left and right limit at zero, for every $\delta>0$
we can choose $\eta<1$ small enough so that $|g(x)-c^+|<\delta$
for $0<x\leq\eta$ and $|g(x)-c^-|<\delta$ for $-\eta\leq x <0$.
Then, using the explicit form of $h$,
\begin{align*}
\lim_{\varepsilon\downarrow 0} \gamma^\varepsilon &\leq \lim_{\varepsilon\downarrow 0}\left\{\int_\varepsilon^\eta \frac{(c^+-\delta)(\varepsilon-x)dx}{|x|^{1+\alpha}} +\int_{-\eta}^{-\varepsilon} \frac{(c^++\delta)(-\varepsilon-x)dx}{|x|^{1+\alpha}} \right\}\\
\lim_{\varepsilon\downarrow 0} \gamma^\varepsilon &\geq
\lim_{\varepsilon\downarrow 0}\left\{\int_\varepsilon^\eta
\frac{(c^++\delta)(\varepsilon-x)dx}{|x|^{1+\alpha}}
+\int_{-\eta}^{-\varepsilon}
\frac{(c^+-\delta)(-\varepsilon-x)dx}{|x|^{1+\alpha}} \right\}
\end{align*}
Explicit evaluation of these integrals together with the fact that
the choice of $\delta$ is arbitrary, yields
$$
\lim_{\varepsilon\downarrow 0} \gamma^\varepsilon =
-\frac{c_+-c_-}{\alpha(\alpha-1)},
$$
and it is easy to check that the third component of the characteristic
triplet of a Lévy process with Lévy density \eqref{stabledens} equals $
-\frac{c_+-c_-}{\alpha(\alpha-1)}$ with the
truncation function $h(x)
= -1 \vee x \wedge 1$  if and only if it equals zero with $h(x)=x$.

\paragraph{Condition \eqref{gammacond} in part 4} We rewrite
 $\gamma^\varepsilon$ as
$$
\gamma^\varepsilon = \gamma + \int_0^\infty \frac{g(x)-g(-x)}{x^2}
\{\varepsilon h(x/\varepsilon)-h(x)\}dx
$$
and apply the dominated convergence, using the fact that
$|\varepsilon h(x/\varepsilon)-h(x)|\leq h(x)$.

\paragraph{Condition \eqref{gammacond} in part 5}
Using the fact that the process has zero drift, we get 
$\gamma^\varepsilon =
\varepsilon^\alpha \int_{\mathbb R} \nu(dx) h(x/\varepsilon)$, and
once again, choosing $\eta>0$ such that $g$ is bounded on
$[-\eta,\eta]$, we get, by dominated convergence:
$$
\lim_{\varepsilon\downarrow 0} \gamma^\varepsilon =
\lim_{\varepsilon\downarrow 0} \varepsilon^\alpha\int_{|x|\leq
\eta}\frac{g(x)h(x/\varepsilon)dx}{|x|^{1+\alpha}} =
\lim_{\varepsilon\downarrow 0}\int_{|x|\leq \eta/\varepsilon}
\frac{g(x\varepsilon)h(x)dx}{|x|^{1+\alpha}} = \int_{\mathbb R}
h(x)\nu^*(dx).
$$

\subsection{Proof of Proposition \ref{cvgtimes.prop}}

\paragraph{Part 1} This will follow if we show that the mapping which to a trajectory $\alpha \in \mathbb D$ (space of càdlàg trajectories) associates $\big(\tau^\alpha_1,\alpha(\tau^\alpha_1)\big)$, with $\tau^\alpha_1:= \inf\{t\geq 0: |\alpha(t)|\geq 1\}$, is continuous in Skorohod topology.
We work component by component. We start with the first component
and study the continuity of the mapping which to a trajectory
$\alpha \in \mathbb D$ associates $\tau^\alpha_1$. This in turn
follows from Proposition VI.2.11 in \cite{jacodshiryaev}, provided
that we prove that the processes $X^\varepsilon$ for every
$\varepsilon$ and $X^*$ satisfy two regularity properties:
\begin{align}
&\inf\{t\geq 0: |Z_t |\geq 1\} = \lim_{\delta\downarrow 1}\inf\{t\geq 0: |Z_t |\geq \delta\}\label{kypr2}\\
&\inf\{t\geq 0: |Z_t |\geq 1\}\leq \inf\{t\geq 0: |Z_{t-} |\geq
1\} \label{kypr1}
\end{align}
almost surely, where $Z$ stands for $X^\varepsilon$ or $X^*$. From
the proof of Lemma 7.10 in \cite{kyprianou}, it follows that
property \eqref{kypr2} holds for every Lévy process unless it is
of compound Poisson type, which is excluded by the conditions of
Proposition \ref{cvginlaw.prop}.\\

\noindent To show Property \eqref{kypr1}, we introduce $\tau =
\inf\{t\geq 0: |Z_{t-} |\geq 1\}$. Remark that $\tau$ is a
stopping time as the hitting time of a Borel set by a càglàd
adapted process (debut theorem). Property \eqref{kypr1} may fail
only if the process $Z$ creeps up to the boundary of $[-1,1]$ and
then immediately jumps back inside this domain, which happens only
if $|Z_{\tau-}|=1$ and $\Delta Z_\tau \neq 0$. Introduce the
sequence $\tau_n = \inf\{t\geq 0: |Z_{t-} |\geq 1-1/n\}$, which
satisfies $\tau_n\leq \tau$. On the set $\{|Z_{\tau-}|=1\}$ also
$\tau_n < \tau$ for all $n$ and it is clear that $\tau_n \to
\tau$. If $|Z_{\tau-}|\neq 1$ it means that the level $1$ is
attained by a jump, and hence $|Z_{\tau-}|< 1$ and $\tau_n=\tau$
as soon as $1-1/n > |Z_{\tau-}|$ so that also $\tau_n \to \tau$.
Therefore, by Proposition I.7 in \cite{bertoin}, on the set
$\{|Z_{\tau-}|=1\}$, $\Delta Z_\tau = 0$.\\

\noindent The continuity of the second component follows from the
proof of Proposition 2.12 in \cite{jacodshiryaev} (part c.)
together with the inequality \eqref{kypr1}.

\paragraph{Part 2} We will show that the family $(\tau^\varepsilon_1)_{\varepsilon > 0}$ has a uniformly bounded exponential moment, which will imply uniform integrability and convergence of $E[(\tau^\varepsilon_1)^kf(X^\varepsilon_{\tau^\varepsilon_1})]$. We treat separately Conditions 1, 2 and 3--5 of Proposition \ref{cvginlaw.prop}.
\paragraph{Condition 1} Since any jump $\Delta X^\varepsilon_t$ with $|\Delta X^\varepsilon|\geq 2$ immediately takes the process $X^\varepsilon$ out of the domain $(-1,1)$, the exit time $\tau^\varepsilon_1$ is dominated by $\tilde \tau^\varepsilon_1 := \inf\{t>0: |\tilde X^\varepsilon_t|\geq 1\}$, where the process $\tilde X^\varepsilon$ is obtained from $X^\varepsilon$ by truncating all jumps greater than $2$ in absolute value. The characteristic exponent of $\tilde X^\varepsilon$ is
$$
\psi_\varepsilon(u) = -\frac{Au^2}{2} + iu\gamma_\varepsilon +
\int_{|x|<2} (e^{iux}-1-iux)\nu^\varepsilon(dx),
$$
where for simplicity we have assumed that the truncation function
satisfies $h(x)=x$ for $|x|<2$. This can be rewritten as
$$
\psi_\varepsilon(u) = -\frac{Au^2}{2} + iu\gamma_\varepsilon +
\int_{|x|<2\varepsilon} \varepsilon^2
(e^{iux/\varepsilon}-1-iux/\varepsilon)\nu(dx) := -\frac{Au^2}{2}
+ \tilde \psi_\varepsilon(u),
$$
and it is easily seen that
$$
|\tilde \psi_\varepsilon(u)|\leq |\gamma_\varepsilon| |u| +
\frac{|u|^2}{2}e^{2|u|} \int_{|x|<2\varepsilon} x^2 \nu(dx),\quad u
\in \mathbb C.
$$
Since $\gamma_\varepsilon \to 0$ as $\varepsilon \to 0$ (see the
proof of Proposition \ref{cvginlaw.prop}), we can find
$\varepsilon_0>0$ such that for all $\varepsilon<\varepsilon_0$
and for all $u\in \mathbb C$ with $|u|=\frac{1}{2}$,
$\frac{2}{A}|\tilde \psi_\varepsilon(u)|<\frac{1}{8}$. From this
bound we deduce: $\Im \psi_\varepsilon(e^{i\pi/12}/2) \leq 0$ and
$\Im \psi_\varepsilon(e^{-i\pi/12}/2) \geq 0$. From the continuity
of $\psi_\varepsilon$ it follows that there exists $\theta \in
\big[-\frac{\pi}{12},\frac{\pi}{12}\big]$ such that
$u^*:=e^{i\theta}/2$ satisfies $\Im \psi_\varepsilon(u^*) = 0$ and
$\Re \psi_\varepsilon(u^*) \in [-\frac{3A}{12},-\frac{A}{16}]$.\\

\noindent Consider now the (complex) exponential martingale
$M^\varepsilon_t = e^{iu^* \tilde X^\varepsilon_t - t
\psi_\varepsilon(u^*)}$. Since $|\tilde
X^\varepsilon_{\tilde\tau^\varepsilon_1 \wedge t}|\leq 3$, we
get that $E[M^\varepsilon_{\tilde\tau^\varepsilon_1}]=1$, and
taking the real part,
$$
E[e^{- \tilde\tau^\varepsilon_1 \psi_\varepsilon(u^*)}] \leq
\frac{e^{3|u^*|}}{\cos(3|u^*|)} = \frac{e^{3/2}}{\cos(3/2)},
$$
which implies
$$
E[e^{\frac{A}{16} \tilde\tau^\varepsilon_1}] \leq
\frac{e^{3/2}}{\cos(3/2)}
$$
for all $\varepsilon<\varepsilon_0$.
\paragraph{Condition 2}
Without loss of generality let $\gamma_0>0$. We use the Lévy-Itô
decomposition of $X$:
$$
X_t = \gamma_0 t + \int_0^t \int_{\mathbb R} z J(ds\times dz),
$$
where $J$ is the jump measure of $X$, and we denote
$$
\tilde X_t := \int_0^t \int_{|z|<2\varepsilon} z J(ds\times dz).
$$
Since any jump $\Delta X$ with $|\Delta X|\geq 2\varepsilon$
immediately takes the process $X^\varepsilon$ out of the domain
$(-1,1)$, for every $k>1$
$$
P\left[\tau^\varepsilon_1> \frac{k}{\gamma_0}\right] \leq
P\left[\gamma_0 t + \tilde X_t \in(-\varepsilon,\varepsilon),
\forall t\leq \frac{k\varepsilon}{\gamma_0}\right] \leq P\left[
|\tilde X_{\frac{k\varepsilon}{\gamma_0}}| > \varepsilon(k-1)
\right].
$$
Since $\tilde X$ has bounded jumps, all its exponential moments
are finite, and therefore for all $\alpha>0,\beta>0$ and $t>0$,
$$
P\left[|\tilde X_t|\geq \alpha\right] \leq e^{-\alpha\beta}E\left[e^{\beta |\tilde
X_t|}\right] \leq e^{-\alpha \beta}
\exp\left(t\int_{|z|<2\varepsilon}(e^{\beta |z|}-1)\nu(dz)\right).
$$
Taking $\alpha = \varepsilon(k-1)$, $\beta =
\frac{1}{\varepsilon}$ and $t=\frac{k\varepsilon}{\gamma_0}$
yields
\begin{align*}
P\left[ |\tilde X_{\frac{k\varepsilon}{\gamma_0}}| >
\varepsilon(k-1) \right] &\leq e^{1-k}
\exp\left(\frac{k\varepsilon}{\gamma_0} \int_{|z|<2\varepsilon}
(e^{|z|/\varepsilon}-1) \nu(dz)\right) \\&\leq e^{1-k}
\exp\left(\frac{k e^2}{\gamma_0} \int_{|z|<2\varepsilon} |z|
\nu(dz)\right).
\end{align*}
Since $X$ is a finite variation process, $\int_{|z|\leq 1}
|z|\nu(dz)<\infty$ and $\lim_{\varepsilon\downarrow 0} \int_{|z|<
2\varepsilon} |z|\nu(dz)$, which means that there exist
$\varepsilon_0>0$, and two constants $c>0$ and $C>0$ such that for
all $\varepsilon\leq \varepsilon_0$ and all $k>1$,
$$
P\left[\tau_\varepsilon> \frac{k}{\gamma_0}\right] \leq C e^{-ck},
$$
which ensures the uniform integrability.

\paragraph{Conditions 3--5}  For $T>0$, the event $\{\tau^\varepsilon_1 >T\}$ occurs only if the process $X^\varepsilon$ does not have any jumps greater or equal to 2 in absolute value on $[0,T]$. Therefore,
$$
P[\tau^\varepsilon_1 > T] \leq \exp\{-T\varepsilon^\alpha
\nu((-\infty,-2\varepsilon]\cup [2\varepsilon,+\infty))\}.
$$
On the other hand,
$$
\varepsilon^\alpha \nu((-\infty,-2\varepsilon]\cup
[2\varepsilon,+\infty)) = \int_{-\infty}^{-2}\frac{g(\varepsilon
x)dx}{|x|^{1+\alpha}} + \int_{2}^{+\infty}\frac{g(\varepsilon
x)dx}{|x|^{1+\alpha}}
$$
is uniformly bounded from below because $g$ has right and left
limits at zero, at least one of which is positive.

\subsection{Proof of Proposition \ref{rate_exit_times}}

\paragraph{Part 1} We first prove the rate of convergence for the first exit time. 
Let $u(x) := E^x[\tau^*_1] = \frac{1-x^2}{2A}1_{|x|\leq 1}$.
Then,
$$
E[\tau^\varepsilon_1 - \tau^*_1] =
E[u(X^\varepsilon_{\tau^\varepsilon_1})+\tau^\varepsilon_1 -
u(0)].
$$
By Itô formula (whose application can be justified, e.g., by
regularizing the function $u$),
\begin{align*}
u(X^\varepsilon_{\tau^\varepsilon_1})+\tau^\varepsilon_1 - u(0) &=
\int_0^{\tau^\varepsilon_1} u'(X^\varepsilon_{t-})
dX^\varepsilon_t + \int_0^{\tau^\varepsilon_1} \left(\frac{A}{2}
u''(X^\varepsilon_t)+1\right) dt \\ &+ \sum_{t\leq
\tau^\varepsilon_1: \Delta X^\varepsilon_t \neq
0}(u(X^\varepsilon_t) -  u(X^\varepsilon_{t-})-\Delta
X^\varepsilon_t  u'(X^\varepsilon_{t-}))\\ &=
\int_0^{\tau^\varepsilon_1} u'(X^\varepsilon_{t-})
dX^\varepsilon_t \\ &+ \sum_{t\leq \tau^\varepsilon_1: \Delta
X^\varepsilon_t \neq 0}(u(X^\varepsilon_t) -
u(X^\varepsilon_{t-})-\Delta X^\varepsilon_t
u'(X^\varepsilon_{t-})),
\end{align*}
where we used the fact that $\frac{A}{2}u''+1=0$. Taking the
expectation, using the boundedness of $u$ and $u'$ and the fact
that the jumps of $X$ have finite variation, we get:
\begin{align*}
E[\tau^\varepsilon_1 - \tau^*_1] = E\left[-\frac{\varepsilon
\gamma_0}{A} \int_0^{\tau^\varepsilon_1} X^\varepsilon_t dt  +
\int_0^{\tau^\varepsilon_1} \int_{\mathbb R} \{u(X^\varepsilon_t +
z)-u(X^\varepsilon_t)\}\nu^\varepsilon(dz)dt\right],
\end{align*}
where $\gamma_0 = \gamma - \int_{\mathbb R}h(x)\nu(dx)$ is the
drift of $X$. Since the limiting process $X^*$ is continuous in
this case, using the Skorokhod representation theorem together with
the fact that the convergence in Skorokhod topology implies convergence
in the local uniform topology (see Theorem VI.1.17 in
\cite{jacodshiryaev}), we get,
$$
\int_0^{\tau^\varepsilon_1} X^\varepsilon_t dt \to
\int_0^{\tau^*_1} X^*_t dt
$$
in law as $\varepsilon\to 0$. Since $\tau^\varepsilon_1$ is uniformly integrable
and $|X^\varepsilon_t|\leq 1$ before $\tau^\varepsilon_1$, also,
$$
\lim_{\varepsilon\downarrow 0} E\left[\int_0^{\tau^\varepsilon_1}
X^\varepsilon_t dt\right] = E\left[\int_0^{\tau^*_1} X^*_t dt\right]=0,
$$
because $X^*$ is a Brownian motion which is a symmetric process.
For the second term under the expectation, we get:
\begin{multline*}
\varepsilon^{-1}E\left[ \int_0^{\tau^\varepsilon_1} \int_{\mathbb
R} \{u(X^\varepsilon_t +
z)-u(X^\varepsilon_t)\}\nu^\varepsilon(dz)dt\right]\\
= E\left[ \int_0^{\tau^\varepsilon_1} \int_{\mathbb R} \varepsilon
\{u(X^\varepsilon_t +
z/\varepsilon)-u(X^\varepsilon_t)\}\nu(dz)dt\right]
\end{multline*} which can be shown to go to zero using the
boundedness of $u$ and $u'$.\\

\noindent To compute the convergence rate of the overshoot, we proceed along the same lines, with the function $u$ now defined by $u(x) = f(x)$ for $|x|\geq 1$ and $u(x) = f(1)$ for $|x|<1$.
  
\paragraph{Part 2} Once again, we start with the first exit time. 
Without loss of generality, assume $\gamma_0>0$. In this case,
$\tau_1^* = \frac{1}{\gamma_0}$. The process $X^\varepsilon$ exits
the interval $(-1,1)$ a.s. in finite time, and we denote by
$U\subset \Omega$ the set of trajectories on which it exits
through the upper barrier. Then,
$$
E[\tau^\varepsilon_1] - E[\tau^*_1] =
E[(\tau^\varepsilon_1-1/\gamma_0)1_{U}] +
E[(\tau^\varepsilon_1-1/\gamma_0)1_{U^c}]
$$
and we analyze the two terms separately. For the first term,
\begin{align*}
\big| E[(\tau^\varepsilon_1-1/\gamma_0)1_{U}]\big| &= \frac{1}{\gamma_0}\big| E\big[\big(X^\varepsilon_{\tau^\varepsilon_1}-1 - \sum_{t\leq \tau^\varepsilon_1} \Delta X^\varepsilon_t\big)1_U\big]\big| \\ &\leq \frac{1}{\gamma_0} E\big[\sum_{t\leq \tau^\varepsilon_1} |\Delta X^\varepsilon_t |\wedge 2\big] = E[\tau^\varepsilon_1] \int_{\mathbb R} (|x|\wedge 2) \nu^\varepsilon(dx)\\
& = E[\tau^\varepsilon_1] \varepsilon \int_{\mathbb R}
(|x/\varepsilon|\wedge 2) \nu(dx),
\end{align*}
where the inequality is due to the fact that on $U$,
$|X^{\varepsilon}_{\tau^\varepsilon_1}-1| \leq \Delta
X^\varepsilon_{\tau^\varepsilon_1}$. Then,
\begin{align*}
&\varepsilon \int_{\mathbb R} (|x/\varepsilon|\wedge 2) \nu(dx) = 2\varepsilon \int_{|x|> 2\varepsilon} \nu(dx) + \int_{|x|\leq 2\varepsilon} |x| \nu(dx)\\
&\qquad \qquad \leq (2\varepsilon)^{1-\beta} \int_{|x|> 2\varepsilon}(|x|^\beta
\wedge 1) \nu(dx) + (2\varepsilon)^{1-\beta}\int_{|x|\leq
2\varepsilon} |x|^\beta \nu(dx),
\end{align*}
from which the result for the first term follows.\\

\noindent To treat the second term, we first estimate the
probability of the set $U^c$. If $\sum_{t\leq 2/\gamma_0} |\Delta
X^\varepsilon_t|\leq 1$ then the process $X^\varepsilon$ surely
exits from the interval $(-1,1)$ through the upper barrier before
time $2/\gamma_0$. Therefore, by the Markov inequality,
\begin{align*}
P[U^c] &\leq P\big[\sum_{t\leq 2/\gamma_0} |\Delta
X^\varepsilon_t|> 1\big] \\ &\leq
P\big[\sum_{t\leq 2/\gamma_0} |\Delta X^\varepsilon_t|1_{|\Delta X^\varepsilon_t|\leq 1}> 1\big] + P\big[\exists t\in [0,2/\gamma_0] : |\Delta X^\varepsilon_t|>1\big] \\
&\leq  E \big[\sum_{t\leq 2/\gamma_0} |\Delta X^\varepsilon_t|1_{|\Delta X^\varepsilon_t|\leq 1} \big] + 1-\exp\big(\frac{2}{\gamma_0} \nu^\varepsilon((-\infty,-1)\cup(1,\infty))\big) \\
&\leq \frac{2}{\gamma_0} \int_{|x|\leq 1} x \nu^{\varepsilon}(dx)
+ \frac{2}{\gamma_0} \int_{|x|> 1}  \nu^{\varepsilon}(dx)  \\ &=
\frac{2}{\gamma_0} \int_{|x|\leq \varepsilon } x \nu(dx) + \frac{2
\varepsilon}{\gamma_0} \int_{|x|> \varepsilon}
\nu(dx) = O(\varepsilon^{1-\beta}).
\end{align*}
The estimate for $E[(\tau^\varepsilon_1-1/\gamma_0)1_{U^c}]$ now
follows by Cauchy-Schwarz inequality and Proposition
\ref{cvgtimes.prop}.\\

\notag We now move to the convergence rate for the overshoot. Let $u(x) = f(x)$ for $|x|\geq 1$ and $u(x) = f(-1) + \frac{x+1}{2}(f(1)-f(-1))$ for $|x|<1$. Applying the Itô formula to $f(X^\varepsilon_{\tau^\varepsilon_1})$ and taking the expectation, we get
\begin{align*}
&E[f(X^\varepsilon_{\tau^\varepsilon_1})-f(X^*_{\tau^*_1})] = E[f(X^\varepsilon_{\tau^\varepsilon_1})-f(1)] \\&= \frac{f(1)-f(-1)}{2} \gamma_0 \{E[\tau^\varepsilon_1] - E[\tau^*_1]\} + E\left[\int_0^{\tau^\varepsilon_1} \int_{\mathbb R} \{f(X^\varepsilon_s+z)-f(X^\varepsilon_s)\}\nu^\varepsilon(dz) ds\right],
\end{align*}
from which the result follows using the boundedness and the Lipschitz property of $f$ and the convergence rate of the first exit time obtained above.

\paragraph{Part 3} Again, we start with the exit time.\\
\noindent\textit{Step 1.} \quad Let $\xi \in (0,1)$ such that
$c<g(x)<C$ for two constants $c$ and $C$ with $0<c<C<\infty$ and
all $x:|x|\leq \xi$. Let $\bar g$ be such that $\bar g(x) = g(x)$
for all $x$ with $|x|\leq \xi$ and $c<g(x)<C$ for all $x$, let
$\bar \nu(dx):= \frac{\bar g(x)}{|x|^{1+\alpha}}dx$ and let $\bar
J$ be a Poisson random measure with intensity $\bar\nu(dx) \times
dt$ independent from $J$. We define the processes $\bar X$ and
$\hat X$ by
\begin{align*}
\bar X_t &:= \big(\gamma - \int_{|x|>\xi} h(x)\nu(dx)\big)t + \int_0^t \int_{|x|\leq \xi} x\tilde J(ds\times dx) + \int_0^t \int_{|x|> \xi} x\bar J(ds\times dx),\\
\hat X_t &:= \big(\gamma - \int_{|x|>\xi} h(x)\nu(dx)\big)t +
\int_0^t \int_{|x|\leq \xi} x\tilde J(ds\times dx).
\end{align*}

\noindent Let $\bar X^\varepsilon_t :=\varepsilon^{-1} \bar
X_{\varepsilon^\alpha t}$, $\hat X^\varepsilon_t
:=\varepsilon^{-1} \hat X_{\varepsilon^\alpha t}$ and let $\bar
\tau^\varepsilon_1$ and $\hat \tau^\varepsilon_1$ be the
corresponding first exit times. By construction, if $\varepsilon \leq
\xi/2$, 
$\tau^\varepsilon_1 \leq \hat \tau^\varepsilon_1$ and if
$\tau^\varepsilon_1 < \hat \tau^\varepsilon_1$ then
$\tau^\varepsilon_1$ is the time of the first jump of $J$ which is
greater than $\xi$ in absolute value; the same statement holds
if $\tau^\varepsilon_1$ is replaced with $\bar\tau^\varepsilon_1$. Let $\mu^\varepsilon$ be the law of
$\hat\tau^\varepsilon_1$. It follows that
\begin{align*}
\big| E[\bar\tau^\varepsilon_1 - \tau^\varepsilon_1]\big| &\leq E[\hat\tau^\varepsilon_1 - \tau^\varepsilon_1]+E[\hat\tau^\varepsilon_1 - \bar\tau^\varepsilon_1]\\
&\leq \int_0^\infty \mu^\varepsilon(dt) t \left(1-e^{-t \varepsilon^\alpha
  \nu(\{x:|x|>\xi\})}\right) \\& \qquad + \int_0^\infty \mu^\varepsilon(dt) t \left(1-e^{-t \varepsilon^\alpha \bar\nu(\{x:|x|>\xi\})}\right)\\
&\leq \varepsilon^\alpha (\nu(\{x:|x|>\xi\})+\bar\nu(\{x:|x|>\xi\}))\int_0^\infty t^2 \mu^\varepsilon(dt)\\
&=\varepsilon^\alpha (\nu(\{x:|x|>\xi\})+\bar
\nu(\{x:|x|>\xi\}))E[(\hat\tau^\varepsilon_1)^2].
\end{align*}
Applying the Proposition \ref{cvgtimes.prop} to the process $\hat
X^\varepsilon$, we get that $E[(\hat\tau^\varepsilon_1)^2]$ is
bounded, and therefore, $\big|E[\bar\tau^\varepsilon_1 -
\tau^\varepsilon_1]\big| = O(\varepsilon^\alpha)$.\\

\noindent\textit{Step 2.} \quad In view of Step 1, it is
sufficient to show that
$$
\lim_{\varepsilon \downarrow
0}\varepsilon^{-\alpha/2}(E[\bar\tau^\varepsilon_1] -
E[\tau^*_1])=0.
$$
Let $\mathbb P^\varepsilon$ be the probability measure under which
the canonical process, denoted by $X$, follows the same law as
$\bar X^\varepsilon$, and $\mathbb P^*$ be the probability measure
under which $X$ follows the same law as $X^*$ By Theorem 33.2 in
\cite{sato}, the restrictions of $\mathbb P^\varepsilon$ and
$\mathbb P^*$ on every finite interval $[0,T]$ are equivalent with
density given by
$$
\frac{d\mathbb P^\varepsilon}{d\mathbb P^*}|_{\mathcal F_T} =
F^\varepsilon_T = \mathcal E(U^\varepsilon)_T,\quad
U^\varepsilon_T = \int_0^T\int_{\mathbb R}
(e^{\phi_\varepsilon(x)}-1)\tilde J^{P^*}(dt\times dx),
$$
where $\tilde J^{P^*}$ is the compensated jump measure of $X$
under $\mathbb P^*$, $\mathcal E$ denotes the Doléans-Dade
exponential, and $\phi^\varepsilon(x):=\frac{\bar g(\varepsilon
x)}{c_+1_{x>0} + c_-1_{x<0}}$.\\

\noindent We denote by $\tau_1$ the first exit time of the
canonical process out of the interval $(-1,1)$ and by
$E^\varepsilon$ and $E^*$ the expectations under the corresponding
probabilities. Let $q\in (1\vee \alpha/\theta,2)$ and $p$ such
that $\frac{1}{q}+\frac{1}{p}=1$. Then by the monotone convergence theorem
and Hölder's inequality,
\begin{align}
|E^\varepsilon[\tau_1] - E^*[\tau_1]| =
|E^*[\tau_1(F^\varepsilon_{\tau_1}-1)]| \leq E^*[(
\tau_1)^p]^{1/p}
E^*[|F^\varepsilon_{\tau_1}-1|^q]^{1/q}\label{holder.eq}
\end{align}
The first factor does not depend on $\varepsilon$ and is clearly
finite ($\tau^*_1$ has an exponential moment). As for the second
factor, since $F^\varepsilon_{t}-1$ is a $\mathbb P^*$-martingale
starting from zero, by the Burkholder-Davis-Gundy inequality we
get,
\begin{align}
&E^*[|F^\varepsilon_{\tau_1}-1|^q] \leq CE^*\big[[F^\varepsilon]^{q/2}_{\tau_1}\big] = CE^*\big[\big(\sum_{t\leq \tau_1:\Delta U^\varepsilon_t\neq 0}(F^\varepsilon_{t-})^2 (\Delta U^\varepsilon_t)^2 \big)^{q/2}\big]\\
&\quad\leq C E^*\big[\sum_{t\leq  \tau_1:\Delta
U^\varepsilon_t\neq 0}(F^\varepsilon_{t-})^q (\Delta
U^\varepsilon_t)^q \big] = C
E^*\big[\int_0^{\tau_1}(F_t^\varepsilon)^q dt\big] \int_{\mathbb
R}(e^{\phi_\varepsilon(x)}-1)^q \nu^*(dx).\label{bdg.eq}
\end{align}
The second factor satisfies
$$
\int_{\mathbb R}\left(e^{\phi_\varepsilon(x)}-1\right)^q \nu^*(dx) =
\varepsilon^\alpha \int_{\mathbb R}\left(e^{\phi_1(x)}-1\right)^q \nu^*(dx) =
O(\varepsilon^\alpha)
$$
by the H\"older property of $g$. For the first factor we get:
\begin{align*}
&E^*\left[\int_0^{\tau_1}(F_t^\varepsilon)^q dt\right] \leq E^*[\tau_1] + E^*\left[\int_0^{ \tau_1}(F_t^\varepsilon)^2 dt\right]\\ &= E^*[\tau_1]  + \int_0^\infty E^*[(F_t^\varepsilon)^2 1_{t\leq  \tau_1}]dt= E^*[\tau_1]  + \int_0^\infty E^\varepsilon[F_t^\varepsilon
1_{t\leq \tau_1}]dt
\end{align*}
To get rid of the stochastic exponential in the last expression,
we would like to make another change of probability measure. Since
$F^\varepsilon$ is not a martingale under $\mathbb P^\varepsilon$,
we represent it as
$$
F^\varepsilon_t = \bar F^\varepsilon_t
\exp\big(tC_\varepsilon\big),
$$
where $\bar F^\varepsilon$ is the Dol\'eans-Dade exponential of
$$
\bar U^\varepsilon_t = \int_0^t\int_{\mathbb R}
(e^{\phi_\varepsilon(x)}-1)\tilde J^{P^\varepsilon}(dt\times dx),
$$
and
\begin{align*}
C_\varepsilon &= \int_{\mathbb R}
\big((e^{\phi_\varepsilon(x)}-1)\phi_\varepsilon(x)-e^{\phi_\varepsilon(x)}+1\big)\nu^*(dx)\\
&= \varepsilon^\alpha \int_{\mathbb R}
\big((e^{\phi_1(x)}-1)\phi_1(x)-e^{\phi_1(x)}+1\big)\nu^*(dx) =
O(\varepsilon^\alpha).
\end{align*}
Then,
\begin{align*}
&E^*\left[\int_0^{ \tau_1}(F_t^\varepsilon)^q dt\right] \leq
E^*[\tau_1]  + \bar E^\varepsilon[e^{\tau_1 C_\varepsilon}],
\end{align*}
where $\bar E^\varepsilon$ denotes the expectation under the
probability $\bar{\mathbb P}^\varepsilon$ such that $\frac{d\bar
{\mathbb P}^\varepsilon}{ d{\mathbb P}^\varepsilon}|_{\mathcal
F_t} = \bar F^\varepsilon_t$. Since $C_\varepsilon\to 0$ and
$\varepsilon\to 0$ and $\tau_1$ has an exponential moment under
$\bar{\mathbb P}^\varepsilon$ (the arguments in the proof of
Proposition \ref{cvgtimes.prop}), we conclude that the first
factor in \eqref{bdg.eq} is finite. Combining this with
\eqref{holder.eq}, the proof is completed.\\

\noindent Let us now turn to the convergence rate for the overshoot. We follow the same steps as above. In step 1, we get, using the boundedness of $f$,
$$
|E[f(\bar X^\varepsilon_{\bar\tau^\varepsilon_1})] - E[f(X^\varepsilon_{\tau^\varepsilon_1})]|\leq C\{P[\tau^\varepsilon_1<\hat \tau^\varepsilon_1]+P[\bar\tau^\varepsilon_1<\hat \tau^\varepsilon_1]\} = O(\varepsilon^\alpha).
$$ 
The rest of the proof is carried out in the same way, with some simplifications due to the boundedness of $f$; for example, the Hölder inequality in \eqref{holder.eq} is not needed. 
\subsection{Proof of Theorem \ref{overshootucp.prop}}

Introduce an auxiliary sequence of times
$(\sigma^\varepsilon_i)_{i\geq 0}$ via $\sigma^\varepsilon_0 = 0$
and $\sigma^\varepsilon_{i+1} = \inf\{t>\sigma^\varepsilon_i :
|X_t - X_{\sigma^\varepsilon_i}|\geq \varepsilon\}$ for $i\geq 1$.
The corresponding counting process is denoted by $M^\varepsilon_t
= \sum_{i\geq 1} 1_{\sigma_i \leq t}$, and it clearly satisfies
$V^\varepsilon(1)_t = M^\varepsilon_{S_t}$ for all $t$. We first
treat the convergence of the process $M^\varepsilon_{t}$.\\

\noindent\textit{Step 1.}\quad Define the process
$$
Z^\varepsilon_t  =
\sum_{i=1}^{[\varepsilon^{-\alpha}t]}(\sigma^{\varepsilon}_{i}-\sigma^{\varepsilon}_{i-1}),
$$
where $[x]$ stands for the integer part of $x$. We first show that
$Z^\varepsilon_t \to t E[\tau^*_1]$ in probability for all $t$.
For every $\Delta>0$,
\begin{align*}
P\big[|Z^\varepsilon_t - t E[\tau^*_1]|>\Delta\big] \leq
P\big[|Z^\varepsilon_t - E[Z^\varepsilon_t]|>\frac{\Delta}{2}\big]
+ 1_{|E[Z^\varepsilon_t]-t E[\tau^*_1]|>\frac{\Delta}{2}}.
\end{align*}
The second term converges to zero because $E[Z^\varepsilon_t] =
[\varepsilon^{-\alpha}t]E[\sigma^\varepsilon_1] =
\varepsilon^\alpha[\varepsilon^{-\alpha}t] \times
\varepsilon^{-\alpha}E[\sigma^\varepsilon_1] \to t E[\tau^*_1]$ by
Proposition \ref{cvgtimes.prop}. For the second term, Chebyshev's
inequality yields:
$$
P\big[|Z^\varepsilon_t - E[Z^\varepsilon_t]|>\frac{\Delta}{2}\big]
\leq \frac{4\text{Var}\,Z_t }{\Delta^2} = \frac{4
[\varepsilon^{-\alpha}t]\text{Var}\,\sigma^\varepsilon_1}{\Delta^2}
\to 0,
$$
because by Proposition \ref{cvgtimes.prop},
$\varepsilon^{-2\alpha}\text{Var}\,\sigma^\varepsilon_1\to
\text{Var}\,\tau^*_1$ as $\varepsilon\to 0$.\\

\noindent\textit{Step 2.}\quad We next show that the convergence
takes place uniformly on compact sets in $t$. Recall
first Dini's theorem which states that a deterministic sequence of
nonnegative increasing functions on $\mathbb{R}^{+}$ converging
pointwise to a continuous function also converges locally
uniformly. Now we use the fact that proving convergence in
probability is equivalent to prove that from any subsequence, one
can extract another subsequence converging almost surely. This
together with Dini's theorem and the pointwise convergence in Step
$1$ gives

$$Z^\varepsilon_t\overset{ucp}{\to} t E[\tau^*_1],\text{ as }\varepsilon\rightarrow 0.$$


\noindent\textit{Step 3.}\quad Our next objective is to deduce the
ucp convergence of $M$ from that of $Z$. Let $\Delta>0, T>0$ and $\bar
M>T/E[\tau^*_1]$. Since $Z^\varepsilon_{M^\varepsilon_t
\varepsilon^\alpha}\leq t$ and $Z^\varepsilon_{(1+M^\varepsilon_t)
\varepsilon^\alpha}> t$, we have
$$P[\sup_{t\leq T} |\varepsilon^\alpha M^\varepsilon_t
E[\tau^*_1]-t|>\Delta]$$ is smaller than
$$P[\sup_{t\leq T} \{\varepsilon^\alpha M^\varepsilon_t
E[\tau^*_1]-Z^\varepsilon_{M^\varepsilon_t
\varepsilon^\alpha}\}>\Delta]+P[\sup_{t\leq T}
\{Z^\varepsilon_{(1+M^\varepsilon_t)
\varepsilon^\alpha}-\varepsilon^\alpha
M^\varepsilon_tE[\tau^*_1]\}>\Delta].$$ Thus, for $\varepsilon$
small enough, there exists some $c>0$ such that this is also
smaller than \begin{align*} 2P[M^\varepsilon_T>\bar M
\varepsilon^{-\alpha}] &+2P\big[\sup_{s\leq \bar
M+c}|Z^\varepsilon_{s} -s E[\tau^*_1]|>\Delta -
E[\tau^*_1]\varepsilon^\alpha,M^\varepsilon_T\leq\bar M
\varepsilon^{-\alpha}\big]\\
\leq 2P[Z^\varepsilon_{\bar M}\leq T] &+2P\big[\sup_{s\leq \bar
M+c}|Z^\varepsilon_{s} -s E[\tau^*_1]|>\Delta/2].
\end{align*}
Since $\bar M>T/E[\tau^*_1]$, the convergence of
$Z^\varepsilon_{\bar M}$ to $\bar M E[\tau^*_1]$ implies that
$P[Z^\varepsilon_{\bar M}\leq T]$ goes to zero. This together with
the ucp convergence of $Z^\varepsilon_t$ in Step 3 gives
$$\varepsilon^\alpha M^\varepsilon_t
E[\tau^*_1]\overset{ucp}{\to} t, \text{ as }\varepsilon\rightarrow
0.$$

\noindent\textit{Step 4.}\quad Define the process
$$
\tilde{Z}^\varepsilon_t(f)
=\varepsilon^{\alpha}\sum_{i=1}^{[\varepsilon^{-\alpha}t]}f\big(\varepsilon^{-1}(X_{\sigma^{\varepsilon}_{i}}-X_{\sigma^{\varepsilon}_{i-1}})\big).$$As
in Step 1, we easily show using Proposition \ref{cvgtimes.prop}
that for $t>0$,
$$\tilde{Z}^\varepsilon_t(f)\rightarrow tE[f(X^*_{\tau^*_1})],$$ in
probability.\\

\noindent\textit{Step 5.} Following Step 2, we obtain
$$\tilde{Z}^\varepsilon_t(f)\overset{ucp}{\to} tE[f(X^*_{\tau^*_1})],\text{ as }\varepsilon\rightarrow 0$$
applying Dini's theorem separately for the positive and negative
parts of $f$.\\

\noindent\textit{Step 6.}  Let $\Delta>0$ and $\eta>0$. Since
$\varepsilon^\alpha M^\varepsilon_t E[\tau^*_1]$ tends ucp to $t$,
for big enough $\varepsilon$,
$$P[\sup_{t\leq T}|\varepsilon^{\alpha}M^\varepsilon_tE[\tau^*_1]E[f(X^*_{\tau^*_1})]-tE[f(X^*_{\tau^*_1})]|>\Delta/2]\leq \eta.$$
Thus,
$$P[\sup_{t\leq T}|\tilde{Z}^\varepsilon_{\varepsilon^{\alpha}M^\varepsilon_t}(f)
E[\tau^*_1]-tE[f(X^*_{\tau^*_1})]|>\Delta]$$ is smaller than
$$P[\sup_{t\leq
T}|\tilde{Z}^\varepsilon_{\varepsilon^{\alpha}M^\varepsilon_t}(f)
E[\tau^*_1]-\varepsilon^{\alpha}M^\varepsilon_tE[\tau^*_1]E[f(X^*_{\tau^*_1})]|>\Delta/2]+\eta.$$
Following the same lines as in Step 3, we eventually obtain
$$\tilde{Z}^\varepsilon_{\varepsilon^{\alpha}M^\varepsilon_t}(f)\overset{ucp}{\to} m(f)t,\text{ as }\varepsilon\rightarrow 0.$$

\noindent\textit{Step 7.} Finally we write
\begin{align*}
&P\big[\sup_{t\leq T}|\varepsilon^\alpha V^{\varepsilon}(f)_t-m(f)S_t|>\delta\big] = P\big[\sup_{t\leq T}|\tilde{Z}^\varepsilon_{\varepsilon^{\alpha}M^\varepsilon_{S_t}}(f)-m(f)S_t|>\delta\big]\\
&\leq P\big[\sup_{t\leq
T}|\tilde{Z}^\varepsilon_{\varepsilon^{\alpha}M^\varepsilon_{S_t}}(f)-m(f)S_t|>\delta,
S_T \leq T^*\big] + P[S_T>T^*].
\end{align*}
Choosing first $T^*$ large enough to make $P[S_T>T^*]$ small, we
can then take $\varepsilon$ small enough to make the first term
small as well. This completes the proof.

\subsection{Proof of Theorem \ref{clt}}
In this proof, we assume without loss of generality that the
function $f_1$ is constant such that $f_1(x)=1$.\\

\noindent\textit{Step 1.}\quad Let $\bar R^\varepsilon_t=(\bar
R^\varepsilon_{t,1},\ldots,\bar R^\varepsilon_{t,d})$ be defined
by
$$
\bar R^\varepsilon_{t,j}=
\varepsilon^{-\alpha/2}\big(\tilde{Z}^{\varepsilon}_{\varepsilon^\alpha
M^\varepsilon_t}(f_j)-tm(f_j)\big).
$$
It is in fact sufficient to show that $\bar R^\varepsilon$ tends
to $B$. Indeed, in that case, the sequence $(\bar R^\varepsilon,
S)$ is $C-$tight (see Corollary VI.3.33 in \cite{jacodshiryaev}).
Using the independence of $S$, we obtain the convergence of finite
dimensional law and finally the convergence in law of $(\bar
R^\varepsilon, S)$ to $(B,S)$. Now using Skorohod representation
theorem, we can place ourselves on the probability space on which
this convergence holds almost surely in Skorohod topology. We
conclude using the fact that for $x$ in the $d$ dimensional
Skorohod space and $y$ an increasing function the $1$ dimensional
Skorohod space function, the application $(x,y)\rightarrow (x\circ
y)$ is
continuous at continuous $(x,y)$ in Skorohod topology.\\

\noindent\textit{Step 2.}\quad In this step we study the
convergence of the process
$L^\varepsilon_t=(L^\varepsilon_{t,1},\ldots,L^\varepsilon_{t,d})$
defined by
$$
L^\varepsilon_{t,j}=\varepsilon^{-\alpha/2}\big(\tilde{Z}^{\varepsilon}_{t/E[\tau^*_1]}(f_j)-m(f_j)Z^{\varepsilon}_{t/E[\tau^*_1]}\big).
$$
We write
$$L^\varepsilon_{t,j}=\sum_{i=1}^{[t/(E[\tau^*_1]\varepsilon^{\alpha})]}\xi_{i,j}^{\varepsilon},$$
 with
$$
\xi_{i,j}^{\varepsilon}=\varepsilon^{\alpha/2}f_j\big(\varepsilon^{-1}(X_{\sigma^{\varepsilon}_{i}}-X_{\sigma^{\varepsilon}_{i-1}})\big)-\varepsilon^{-\alpha/2}m(f_j)(\sigma^{\varepsilon}_{i}-\sigma^{\varepsilon}_{i-1}).
$$
Using that
$$\big\{\varepsilon^{-1}(X_{\sigma^{\varepsilon}_{i}}-X_{\sigma^{\varepsilon}_{i-1}}),\sigma^{\varepsilon}_{i}-\sigma^{\varepsilon}_{i-1}\big\}$$
and
$\{X^{\varepsilon}_{\tau_1^{\varepsilon}},\varepsilon^{\alpha}\tau^{\varepsilon}_1\}$
have the same law, we get
\begin{align*}
E[\xi_{i,j}^{\varepsilon}]&=\varepsilon^{\alpha/2}\big(E[f_j(X^{\varepsilon}_{\tau_1^{\varepsilon}})]-m(f_j)E[\tau^{\varepsilon}_1]\big)\\
&=\varepsilon^{\alpha/2}\big(E[f_j(X^{\varepsilon}_{\tau_1^{\varepsilon}})]-E[f_j(X^{*}_{\tau_1^{*}})]+m(f_j)(E[\tau^*_1]-E[\tau_1^{\varepsilon}])\big)
\end{align*}
and for $1\leq j,k\leq d$,
$$
E[\xi_{i,j}^{\varepsilon}\xi_{i,k}^{\varepsilon}]=\varepsilon^{\alpha}E\big[\big(f_j(X^{\varepsilon}_{\tau_1^{\varepsilon}})-m(f_j)\tau_1^{\varepsilon}\big)\big(f_k(X^{\varepsilon}_{\tau_1^{\varepsilon}})-m(f_k)\tau_1^{\varepsilon}\big)\big].
$$
Moreover, for some positive constant $c$,
$$E[(\xi_{i,j}^{\varepsilon})^4]\leq c\varepsilon^{2\alpha}.$$
From the specific assumptions on $X$ for Theorem \ref{clt}, we get
$$\sum_{i=1}^{[t/(E[\tau^*_1]\varepsilon^{\alpha})]}E[\xi_{i,j}^{\varepsilon}]\rightarrow
0.$$ Now, using Proposition \ref{cvgtimes.prop}, we obtain
$$\sum_{i=1}^{[t/(E[\tau^*_1]\varepsilon^{\alpha})]}\big(E[\xi_{i,j}^{\varepsilon}\xi_{i,k}^{\varepsilon}]-E[\xi_{i,j}^{\varepsilon}]E[\xi_{i,k}^{\varepsilon}]\big)\rightarrow (t/E[\tau^*_1])C_{j,k}$$ with
$$C_{j,k}=\text{Cov}[f_j(X^{*}_{\tau_1^{*}})-m(f_j)\tau^*_1,f_k(X^{*}_{\tau_1^{*}})-m(f_k)\tau^*_1].$$

\noindent Using a usual theorem on the convergence of
triangular arrays, see Theorem VIII.3.32 in \cite{jacodshiryaev},
we obtain that $L^\varepsilon$ converges in law to a continuous
centered $\mathbb{R}^d-$valued Gaussian process with independent
increments $B$ such that
$E[B_{t,j}B_{t,k}]=(t/(E[\tau^*_1])C_{j,k}$.\\

\noindent\textit{Step 3.}\quad We introduce two families of time
changes converging ucp to identity: $\eta^\varepsilon_t =
\varepsilon^\alpha M^\varepsilon_t E[\tau^*_1]$ and
$\bar\eta^\varepsilon_t = \varepsilon^\alpha (1+M^\varepsilon_t)
E[\tau^*_1]$. Since the ucp convergence implies the convergence in
law in the Skorohod space, the sequences $\eta^\varepsilon_t$ and
$\bar\eta^\varepsilon_t$ are C-tight. The sequence
$L^{\varepsilon}_t$ being also C-tight, the sequence of
$d+2$-dimensional processes
$(L^{\varepsilon}_t,\eta^\varepsilon_t,\bar\eta^\varepsilon_t)$ is
C-tight. Since the time changes converge to deterministic limits,
we also get the finite dimensional convergence of the preceding
sequence which implies its convergence in law in the Skorohod space for the Skorohod topology.\\

\noindent By the Skorohod representation theorem, we can place
ourselves on the probability space on which $L^\varepsilon\to B$,
$\eta^\varepsilon_t\rightarrow t$ and
$\bar\eta^\varepsilon_t\rightarrow t$ almost surely in Skorohod
topology. Using again the continuity of composition by time change
at continuous limits, we get that
$L^{\varepsilon}_{\eta^\varepsilon_t}\rightarrow B_t$ and
$L^{\varepsilon}_{\bar\eta^\varepsilon_t}\rightarrow B_t$. Since
$B_t$ is continuous, this implies
$L^{\varepsilon}_{\eta^\varepsilon_t}-L^{\varepsilon}_{\bar\eta^\varepsilon_t}\rightarrow
0$ and so, using that $f_1(x)=1$,
$$\varepsilon^{{\alpha/2}}+\varepsilon^{-{\alpha/2}}\big(m(f_1)(Z^{\varepsilon}_{\varepsilon^\alpha (M^\varepsilon_t+1)}-Z^{\varepsilon}_{\varepsilon^\alpha M^\varepsilon_t})\big)\rightarrow 0,$$
which gives
$$\varepsilon^{-{\alpha/2}}(Z^{\varepsilon}_{\varepsilon^\alpha M^\varepsilon_t}-Z^{\varepsilon}_{\varepsilon^\alpha (M^\varepsilon_t+1)})\rightarrow 0.$$
This also implies the convergence for the local uniform topology
(see Theorem VI.1.17 in \cite{jacodshiryaev}). Since by
construction $Z^\varepsilon_{M^\varepsilon_t
\varepsilon^\alpha}\leq t$ and $Z^\varepsilon_{(1+M^\varepsilon_t)
\varepsilon^\alpha}> t$ we get
$$|Z^{\varepsilon}_{\varepsilon^\alpha M^\varepsilon_t}-t|\leq|Z^{\varepsilon}_{\varepsilon^\alpha
M^\varepsilon_t}-Z^{\varepsilon}_{\varepsilon^\alpha
(M^\varepsilon_t+1)}|.$$ Thus,
$$\varepsilon^{-{\alpha/2}}(Z^{\varepsilon}_{\varepsilon^\alpha M^\varepsilon_t}-t)\rightarrow 0.$$ Eventually, we use that
$\bar
R^\varepsilon_t=L^{\varepsilon}_{\eta^\varepsilon_t}+\gamma^{\varepsilon}_t,$
with
$$\gamma^{\varepsilon}_{j,t}=m(f_j)\varepsilon^{-{\alpha/2}}(Z^{\varepsilon}_{\varepsilon^\alpha
M^\varepsilon_t}-t).$$ Since  $\gamma^{\varepsilon}_t\rightarrow
0$, the result follows.

\subsection{Proof of Proposition \ref{noclt}}
The idea is to repeat the proof of Theorem
\ref{overshootucp.prop}, using sharper estimates \eqref{beta1} et
\eqref{beta2} to obtain the convergence rate. We only give the sketch
of the proof. \\

\noindent 
Define the process
$$
U^\varepsilon_t = \varepsilon^{-(1-\delta-\beta)\vee -\frac{1}{2}}
\left\{ \sum_{i=1}^{[\varepsilon^{-1}t]} (\sigma^\varepsilon_i -
  \sigma^{\varepsilon}_{i-1}) - tE[\tau^*_1]
\right\}.
$$
We recall that in our setting $\tau^*_1$ is deterministic, but we
stick to the notation of the proof of Theorem
\ref{overshootucp.prop}. Then,
\begin{align*}
U^\varepsilon_t = \varepsilon^{-(1-\delta-\beta)\vee -\frac{1}{2}}
\left\{ \sum_{i=1}^{[\varepsilon^{-1}t]} (\sigma^\varepsilon_i -
  \sigma^{\varepsilon}_{i-1}) - [\varepsilon^{-1} t]
  E[\sigma^\varepsilon_1]\right\} \\ + \varepsilon^{-(1-\delta-\beta)\vee -\frac{1}{2}}\left\{\varepsilon[\varepsilon^{-1} t]
  E[\tau^\varepsilon_1]- tE[\tau^*_1]
\right\}.
\end{align*}
The bound \eqref{beta1} implies that the terms in the second line
converge to zero uniformly in $t$ on compacts. The terms in the first
line, by Kolmogorov's inequality, satisfy
\begin{align*}
&P\left[\sup_{t\leq t_0}\varepsilon^{-(1-\delta-\beta)\vee -\frac{1}{2}}
\left| \sum_{i=1}^{[\varepsilon^{-1}t]} (\sigma^\varepsilon_i -
  \sigma^{\varepsilon}_{i-1}) - [\varepsilon^{-1} t]
  E[\sigma^\varepsilon_1]\right|\geq \lambda\right] \\ &\leq
\frac{1}{\lambda^2} \varepsilon^{-2(1-\delta-\beta)\vee
  -1}[\varepsilon^{-1}t_0]\text{Var}\, \sigma^\varepsilon_1 \leq \frac{t_0}{\lambda^2} \text{Var}\, \tau^\varepsilon_1,
\end{align*}
which converges to zero as $\varepsilon\to 0$ because
$E[\tau^\varepsilon_1] \to E[\tau^*_1]$, $E[(\tau^\varepsilon_1)^2]
\to E[(\tau^*_1)^2]$ and $\tau^*_1$ is deterministic. We have
therefore shown that $U^\varepsilon \xrightarrow{ucp} 0$. \\

\noindent Now we repeat the arguments of step 3 of the proof of
Theorem \ref{overshootucp.prop} to show that 
$$
\varepsilon^{-(1-\delta-\beta)\vee -\frac{1}{2}}\{\varepsilon
M^\varepsilon_t E[\tau^*_1] - t\}\xrightarrow{ucp} 0. 
$$
Finally, we define 
$$
\tilde U^\varepsilon_t(f) = \varepsilon^{-(1-\delta-\beta)\vee -\frac{1}{2}}
\left\{ \sum_{i=1}^{[\varepsilon^{-1}t]} (f(\varepsilon^{-1}(X_{\sigma^\varepsilon_i} -
  X_{\sigma^{\varepsilon}_{i-1}})) - tE[f(X^*_{\tau^*_1})]
\right\}.
$$
and show that $\tilde U^\varepsilon (f)\xrightarrow{ucp} 0$ using the
same argument as above. The proof can then be completed by repeating the steps 5--7 of the proof of
Theorem \ref{overshootucp.prop} with the process $\tilde
Z^\varepsilon(f)$ replaced by $\tilde U^\varepsilon(f)$.


\end{document}